\documentclass{eptcs}

\usepackage{amsthm}
\usepackage{mathtools}

\usepackage[inline]{enumitem}
\usepackage{ifthen}
\usepackage[utf8]{inputenc} 
\usepackage{xcolor}

\usepackage[capitalize]{cleveref}
\usepackage{url}

\usepackage{tikz}

\usepackage{amssymb}
\usepackage[varg,bigdelims]{newpxmath}
\usepackage{mathrsfs}
\usepackage{dutchcal}
\usepackage{ebproof}
\usepackage{stmaryrd}
\usepackage{float}

  \crefformat{enumi}{\card#2#1#3}
  \crefalias{section}{section}

  \hypersetup{final}

  \setlist{nosep}
  \setlistdepth{6}

  \usetikzlibrary{ 
  	cd,
  	math,
  	decorations.markings,
		decorations.pathreplacing,
  	positioning,
  	arrows.meta,
  	shapes,
		shadows,
		shadings,
  	calc,
  	fit,
  	quotes,
  	intersections,
    circuits,
    circuits.ee.IEC
  }
  
  \tikzset{
biml/.tip={Glyph[glyph math command=triangleleft, glyph length=.95ex]},
bimr/.tip={Glyph[glyph math command=triangleright, glyph length=.95ex]},
}

\tikzset{
	tick/.style={postaction={
  	decorate,
    decoration={markings, mark=at position 0.5 with
    	{\draw[-] (0,.4ex) -- (0,-.4ex);}}}
  }
} 
\tikzset{
	slash/.style={postaction={
  	decorate,
    decoration={markings, mark=at position 0.5 with
    	{\node[font=\footnotesize] {\rotatebox{90}{$\approx$}};}}}
  }
}

\newcommand{\xslashar}[1]{\begin{tikzcd}[baseline=-0.5ex,cramped,sep=small,ampersand 
replacement=\&]{}\ar[r,slash, "{#1}"]\&{}\end{tikzcd}}

\tikzset{
	oriented WD/.style={
		every to/.style={out=0,in=180,draw},
    label/.style={
    	font=\everymath\expandafter{\the\everymath\scriptstyle},
      inner sep=0pt,
      node distance=2pt and -2pt},
    semithick,
    node distance=1 and 1,
    decoration={markings, mark=at position \stringdecpos with \stringdec},
    ar/.style={postaction={decorate}},
    execute at begin picture={\tikzset{
    	x=\bbx, y=\bby,
      every fit/.style={inner xsep=\bbx, inner ysep=\bby}}}
    },
    string decoration/.store in=\stringdec,
    string decoration={\arrow{stealth};},
    string decoration pos/.store in=\stringdecpos,
    string decoration pos=.7,
    bbx/.store in=\bbx,
    bbx = 1.5cm,
    bby/.store in=\bby,
    bby = 1.5ex,
    bb port sep/.store in=\bbportsep,
    bb port sep=1.5,
    bb port length/.store in=\bbportlen,
    bb port length=4pt,
    bb penetrate/.store in=\bbpenetrate,
    bb penetrate=0,
    bb min width/.store in=\bbminwidth,
    bb min width=1cm,
    bb rounded corners/.store in=\bbcorners,
    bb rounded corners=2pt,
    bb spider/.style={
    	bb port sep=1, bb port length=10pt, bbx=.4cm, bb min width=.4cm, bby=.8ex},
    bb small/.style={
    	bb port sep=1, bb port length=2.5pt, bbx=.4cm, bb min width=.4cm, bby=.7ex},
		bb medium/.style={
			bb port sep=1, bb port length=2.5pt, bbx=.4cm, bb min width=.4cm, bby=.9ex},
    bb/.code 2 args={
    	\pgfmathsetlengthmacro{\bbheight}{\bbportsep * (max(#1,#2)+1) * \bby}
      \pgfkeysalso{draw,minimum height=\bbheight,minimum
       width=\bbminwidth,outer sep=0pt,
         rounded corners=\bbcorners,thick,
         prefix after command={\pgfextra{\let\fixname\tikzlastnode}},
         append after command={\pgfextra{\draw
            \ifnum #1=0{} \else foreach \i in {1,...,#1} {
            	($(\fixname.north west)!{\i/(#1+1)}!(\fixname.south west)$) +(-\bbportlen,0) coordinate (\fixname_in\i) -- +(\bbpenetrate,0) coordinate (\fixname_in\i')}\fi 
            \ifnum #2=0{} \else foreach \i in {1,...,#2} {
            	($(\fixname.north east)!{\i/(#2+1)}!(\fixname.south east)$) +(-
\bbpenetrate,0) coordinate (\fixname_out\i') -- +(\bbportlen,0) coordinate (\fixname_out\i)}\fi;
           }}}
		},
			bb name/.style={
     	append after command={
				\pgfextra{\node[anchor=north] at (\fixname.north) {#1};}
			}
		}
  }

\tikzset{Yonepart/.pic={
	\node[bb={1}{2},bb name = {\tiny$X_{11}$}] (X11) {};
	\node[bb={2}{2},below right=of X11,bb name = {\tiny$X_{12}$}] (X12) {};
	\node[bb={2}{1}, above right=of X12,bb name = {\tiny$X_{13}$}] (X13) {};
	\node[bb={0}{0}, inner xsep=10pt, fit={($(X11.north west)+(.3,1.5)$) (X12)  ($(X13.east)+(-.3,0)$)},bb name = {\scriptsize $Y_1$}] (Y1) {};
	\coordinate (Y1_in1') at (X11_in1-|Y1.west);
	\coordinate (Y1_in1) at (X11_in1-|Y1.west);
	\coordinate (Y1_in2') at (X12_in2-|Y1.west);
	\coordinate (Y1_in2) at (X12_in2-|Y1.west);
	\coordinate (Y1_out1') at (X13_out1-|Y1.east);
	\coordinate (Y1_out1) at (X13_out1-|Y1.east);
	\coordinate (Y1_out2') at (X12_out2-|Y1.east);
	\coordinate (Y1_out2) at (X12_out2-|Y1.east);
	\draw (Y1_in1') to (X11_in1);	
	\draw (Y1_in2') to (X12_in2);
	\draw (X11_out1) to (X13_in1);
	\draw (X11_out2) to (X12_in1);
	\draw (X12_out1) to (X13_in2);
	\draw (X12_out2) to (Y1_out2');
	\draw (X13_out1) to (Y1_out1');
	\coordinate (bottombox) at ($(X12.south)$);
	\coordinate (rightbox) at ($(X13.east)$);
	\coordinate (Y1northwest) at ($(Y1.north west)$);
	}
}
  \tikzset{Ytwopart/.pic={
	\node[bb={2}{2}, bb name = {\tiny$X_{21}$}] (X21) {};
	\node[bb={1}{2},above right=-1 and 1 of X21,bb name = {\tiny$X_{22}$}] (X22) {};
	\node[bb={0}{0}, inner xsep=10pt, fit={($(X21.south west)+(-.25,0)$) ($(X22.north east)+(.25,3.5)$)},bb name = {\scriptsize$Y_2$}] (Y2){};
	\coordinate (Y2_in1') at (X21_in2-|Y2.west);
	\coordinate (Y2_in1) at (X21_in2-|Y2.west);
	\coordinate (Y2_out1') at (X22_out2-|Y2.east);
	\coordinate (Y2_out1) at (X22_out2-|Y2.east);
	\coordinate (Y2_out2') at (X21_out2-|Y2.east);
	\coordinate (Y2_out2) at (X21_out2-|Y2.east);	
	\draw (Y2_in1') to (X21_in2);
	\draw (X21_out1) to (X22_in1);
	\draw (X22_out2) to (Y2_out1');
	\draw let \p1=(X22.south east), \p2=($(Y2_out2)$), \n1={\y1-\bby}, \n2=\bbportlen in
	  (X21_out2) to (\x1+\n2,\n1) -- (\x1+\n2,\n1) to (Y2_out2');
	\draw let \p1=(X22.north east), \p2=(X21.north west), \n1={\y1+\bby}, \n2=\bbportlen in
          (X22_out1) to[in=0] (\x1+\n2,\n1) -- (\x2-\n2,\n1) to[out=180] (X21_in1);
          }
}

\tikzset{SmallNestingPic/.pic={
  \path (0,0) pic [purple] {Yonepart};
  \path ($(rightbox)+(4,-5)$) pic [blue!40!black] {Ytwopart};
  
  \node[bb={0}{0}, fit={($(Y1northwest)+(-.5,4)$) ($(Y2.south east)+(1,0)$)}, bb name={\small $Z$}] (Z) {};
  \coordinate[above=\bby] (helper) at (Y2.north west);
	\coordinate (Z_in1') at (Y1_in2-|Z.west);
	\coordinate (Z_in1) at (Y1_in2-|Z.west);
	\coordinate (Z_out1') at (helper-|Z.east);
	\coordinate (Z_out1) at (helper-|Z.east);
	\coordinate (Z_out2') at (Y2_out2-|Z.east);
	\coordinate (Z_out2) at (Y2_out2-|Z.east);	  
  \draw (Z_in1') to (Y1_in2);
  \draw let \p1=(Y2.north west),\p2=(Y2.north east),\n1={\y2+\bby},\n2=\bbportlen in
  (Y1_out1) to (\x1+\n2,\n1)--(\x2+\n2,\n1) to (Z_out1');
  \draw (Y1_out2) to (Y2_in1);
  \draw (Y2_out2) to (Z_out2');
  \draw let \p1=(Y2.north east), \p2=(Y1.north west), \n1={\y2+\bby}, \n2=\bbportlen in
          (Y2_out1) to[in=0] (\x1+\n2,\n1) -- (\x2-\n2,\n1) to[out=180] (Y1_in1);
          }
}

\theoremstyle{definition}
\newtheorem{definitionx}{Definition}[section]
\newtheorem{examplex}[definitionx]{Example}

\theoremstyle{plain}

\newtheorem{theorem}[definitionx]{Theorem}
\newtheorem{proposition}[definitionx]{Proposition}

\newtheorem*{theorem*}{Theorem}
\newtheorem*{proposition*}{Proposition}
\newtheorem*{corollary*}{Corollary}
\newtheorem*{lemma*}{Lemma}
\newtheorem*{warning*}{Warning}

\newenvironment{example}
  {\pushQED{\qed}\examplex}
  {\popQED\endexamplex}

  \newenvironment{definition}
  {\pushQED{\qed}\definitionx}
  {\popQED\enddefinitionx}

\DeclareSymbolFont{stmry}{U}{stmry}{m}{n}
\DeclareMathSymbol\fatsemi\mathop{stmry}{"23}

\DeclareFontFamily{U}{mathx}{\hyphenchar\font45}
\DeclareFontShape{U}{mathx}{m}{n}{
      <5> <6> <7> <8> <9> <10>
      <10.95> <12> <14.4> <17.28> <20.74> <24.88>
      mathx10
      }{}
\DeclareSymbolFont{mathx}{U}{mathx}{m}{n}
\DeclareFontSubstitution{U}{mathx}{m}{n}
\DeclareMathAccent{\widecheck}{0}{mathx}{"71}

\DeclareMathOperator{\Hom}{Hom}

\DeclareMathOperator{\ob}{Ob}

\newcommand{\cat}[1]{\mathcal{#1}}
\newcommand{\Cat}[1]{\textbf{#1}}

\newcommand{\id}{\mathrm{id}}
\newcommand{\then}{\mathbin{\fatsemi}}

\newcommand{\too}{\longrightarrow}

\newcommand{\To}[2][]{\xrightarrow[#1]{#2}}

\newcommand{\from}{\leftarrow}
\newcommand{\fromm}{\longleftarrow}

\newcommand{\From}[1]{\xleftarrow{#1}}

\newcommand{\slashar}{\xslashar{}}
\newcommand{\card}{\,^{\#}}

\newcommand{\tn}[1]{\textnormal{#1}}

\newcommand{\nn}{\mathbb{N}}

\newcommand{\rr}{\mathbb{R}}

\newcommand{\smset}{\Cat{Set}}
\newcommand{\smcat}{\Cat{Cat}}

\newcommand{\act}{\tn{act}}
\newcommand{\upd}{\tn{upd}}

\newcommand{\yon}{\mathcal{y}}
\newcommand{\poly}{\Cat{Poly}}

\newcommand{\tri}{\mathbin{\triangleleft}}

\newcommand{\bet}{\Delta^+}

\newcommand{\biglens}[2]{
     \begin{bmatrix}{\vphantom{f_f^f}#2} \\ {\vphantom{f_f^f}#1} \end{bmatrix}
}
\newcommand{\littlelens}[2]{
     \begin{bsmallmatrix}{\vphantom{f}#2} \\ {\vphantom{f}#1} \end{bsmallmatrix}
}
\newcommand{\lens}[2]{
  \relax\if@display
     \biglens{#1}{#2}
  \else
     \littlelens{#1}{#2}
  \fi
}

\newcommand{\slogan}[1]{\begin{center}\textit{#1}\end{center}}

\newcommand{\qqand}{\qquad\text{and}\qquad}

\newcommand{\coalg}{\tn{-}\Cat{Coalg}}

\newcommand{\org}{{\mathbb{O}\Cat{rg}}}
\renewcommand{\S}{{\Cat{S}}}
\newcommand{\idcoalg}[1]{{\{\id_{#1}\}}}

\allowdisplaybreaks
\begin{document}

\title{Dynamic Operads, Dynamic Categories:\\From Deep Learning to Prediction Markets}

\def\titlerunning{Dynamic Categories, Dynamic Operads}
\author{Brandon T. Shapiro \and David I. Spivak}
\def\authorrunning{B.T.\ Shapiro and D.I.\ Spivak}
\date{\vspace{-.2in}}

\maketitle

\begin{abstract}
Natural organized systems adapt to internal and external pressures and this happens at all levels of the abstraction hierarchy. Wanting to think clearly about this idea motivates our paper, and so the idea is elaborated extensively in the introduction, which should be broadly accessible to a philosophically-interested audience.

In the remaining sections, we turn to more compressed category theory. We define the monoidal double category $\org$ of dynamic organizations, we provide definitions of $\org$-enriched, or \emph{dynamic}, categorical structures---e.g.\ dynamic categories, operads, and monoidal categories---and we show how they instantiate the motivating philosophical ideas. We give two examples of dynamic categorical structures: prediction markets as a dynamic operad and deep learning as a dynamic monoidal category.
\end{abstract}

\section{Introduction}

Intuitively, an \emph{open dynamical system} is a machine or worker with an interface by which to interact with whatever else is out there. Open dynamical systems can be organized as circuits or control loops, so that they affect each other by their outward expressions of internal work, and thereby possibly form a more complex worker. The framework here is fractal---or more precisely \emph{operadic}---in its structure: organizations of workers can be nested into arbitrary hierarchies of abstraction.
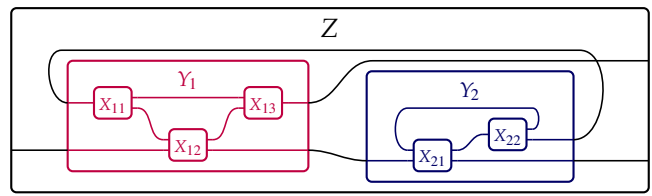
\begin{figure}[H]
\[
\begin{tikzpicture}[oriented WD, bb min width =.5cm, bbx=.5cm, bb port sep =1,bb port length=.08cm, bby=.14cm]
\path (0,0) pic {SmallNestingPic};
\end{tikzpicture}
\]
\caption{
A nesting of interacting open dynamical systems: the $X_{i,j}$ are wired together to form the $Y_i$, which are wired together to form $Z$; typically these groupings are chosen to create new abstractions, e.g.\ in logical circuits or control systems. The permanence of the above-displayed wiring pattern is exactly what is relaxed in this paper; a dynamic organization is one in which interactions may change dynamically based on what flows within the system.
}\label{fig.nesting}
\end{figure}

But if we think about some things that interact to do work in the real world, we notice that often the organization itself---the connections themselves---change. Unlike what we see in \cref{fig.nesting}, the way we connect this hour may be different from the way we connect next hour; in particular, our interfaces go in and out of contact. At the end of this paragraph, look away from the page for a few seconds, think about some things you know that interact together or influence each other, and ask yourself three questions about them: Do these things ever stop interacting? If so, do they ever start interacting again? And how is it decided? 

\subsection{Accounting for organizational change}

We propose that the metaphysical nature and scope of these questions should be complemented by some sort of guard rails to keep our contemplation on track. This is the role of mathematics in our work. It provides a symbolic \emph{accounting system} which is articulate enough to facilitate one person in explicating an example and asking questions about it.

The category $\poly$ of polynomial functors in one variable is an ergonomic mathematical structure with many applications and spin-off categorical gadgets. We will begin in \cref{chap.org} by recalling one such gadget from \cite{spivak2021learners}: a category-enriched multicategory $\org$ that will be the conceptual centerpiece of our accounting system. Its objects are polynomial functors in one variable, and its morphisms are polynomial coalgebras related to a certain monoidal closed structure on $\poly$. We will see that the morphisms in $\org$ are intuitively ``collective organizational patterns that change dynamically''.

Leaving the mathematics aside until \cref{chap.org}---at which point we will have almost nothing more to say about the background philosophy---let's return to the question ``how is the organizational pattern between various systems decided, moment-by-moment?'' Let's mesh this question with the idea that the so-organized systems can be nested into arbitrary hierarchies of abstraction. And let's think about all this in the frame of a certain worldview which we invite you the reader to engage with like a fictional movie, not intended to convince you of fact but instead simply to convey an experience. Here goes.

In this worldview, we notice that everything that makes any sense to us happens to be a collective. A cell body, a human body, an antibody, Topos Institute, an idea, an airport, a sentence, a mathematical definition, a grain of sand, ... each is a collective of interacting parts that may themselves be collectives. 

It's quite often the case that these collectives, like the ship of Theseus, are not permanent organizations that are fixed for all time; they are adapting to forces from within and without the system. Even a grain of sand can break or melt; even a mathematical definition can be refactored. So then what's outside the system, generating these forces that influence it? We imagine that what's outside is in fact more of the same kind of stuff as what's inside, just not as cohesive perhaps. Let's go full-on woo: if the universe is a big system, then maybe the sort of thing that happens in our head is---in some way---just like what happens outside of it. Maybe the motives that organize Brandon and David into a collaborative thinking and paper-writing unit are, in the some reasonable account, of the same nature as the motives that organize each one of them into a body. 

But is this right? How could you check such a claim? One would need to give a reasonable account of it, and since we as authors can't currently give such an account, we don't make this claim. Instead, what we present here is an \emph{accounting system} in which the woo-person, (or would it instead be the reductive materialist?) who thought that what went on inside the head was somehow the same as what went on outside, could endeavor to provide such an account of their thinking.

\subsection{Dynamic categorical structures}

Our main definition in this paper is what we call an \emph{dynamic categorical structure}. We might poetically say that a \emph{dynamic} category is one where the morphisms between two objects change in response to what flows between those objects. To define it, we first refactor the definition of $\org$ from \cite{spivak2021learners} from an operad to a monoidal double category; we then define a dynamic *thing* to be a *thing* enriched in $\org$. Once these are defined, we give a couple examples: a \emph{prediction market} operad and a \emph{deep learning} monoidal category. In the prediction market, a population $Y$ predicts a distribution based on the predictions of its member populations $X_i$ weighted by their reputations, and the reputations change dynamically based on the returned outcome. A similar story holds with deep learning.

We thank you the reader for having postponed your counterpoints and counterexamples, and we ask you to reengage both skepticism and interest as you see fit. We invite you to ask openly: what's not a collective of interacting parts that are themselves collectives? Nature, love, or experience perhaps? It all depends on how you look. What we present here is an accounting system for making sense of a certain sort of experiential pattern; the matter itself is whatever it is.

\subsection{Acknowledgments}

The influences on this paper are too numerous and unranked to name, but in particular we thank Sophie Libkind for stimulating conversations, and we thank Scott Garrabrant for teaching us about Kelly betting, which partially inspired \cref{sec.kelley}. Thanks also to David Jaz Myers and Samantha Jarvis for catching typos.

This material is based upon work supported by the Air Force Office of Scientific Research under award number FA9550-20-1-0348.

\section{The Monoidal Double Category $\org$}\label{chap.org}

In \cite{spivak2021learners}, the second author defined a category-enriched multicategory $\org$, whose objects are polynomials and whose morphisms are polynomial coalgebras. In this section, we describe how $\org$ in fact more naturally takes the form of a monoidal double category, with coalgebras as horizontal morphisms, maps of polynomials as vertical morphisms, and the Dirichlet tensor product $\otimes$ (see \eqref{eqn.tensor} below) providing the monoidal structure.%
\footnote{In fact, $\org$ is a duoidal double category, with a second monoidal structure $\tri$, but we will not use that here.}

Before we begin, recall that a polynomial is a functor $p\colon \smset\to\smset$ which is isomorphic to a sum of representables; following \cite{spivak2021learners}, we denote $p,q\in\poly$ by
\begin{equation}\label{eqn.poly_notation}
p = \sum_{I \in p(1)} \yon^{p[I]} \qqand q = \sum_{J \in q(1)} \yon^{q[J]}
\end{equation}
and refer to each $I\in p(1)$ as a \emph{$p$-position} and to each $i\in p[I]$ as a \emph{$p$-direction at $I$}. A map $\phi\colon p\to q$ of polynomials is a natural transformation. Combinatorially, $\phi$ provides: for each $I\in p(1)$ a choice of $\phi(I)\in q(1)$ and for each $j\in q[\phi(I)]$ a choice of $\phi(I,j)\in p[I]$.%
\footnote{In \cite{spivak2021learners}, what we denote $\phi(I)$ is denoted $\phi_1(I)$ and what we denote $\phi(I,j)$ is denoted $\phi^\sharp_I(j)$.}

For polynomials $p,q$, their Dirichlet tensor product is the polynomial
\begin{equation}\label{eqn.tensor}
p \otimes q = \sum_{(I,J) \in p(1) \times q(1)} \yon^{p[I] \times q[J]}
\end{equation}

\subsection{$[p,q]$-coalgebras}

We first recall the definitions of the internal-hom polynomials $[p,q]$ and concretely describe the category of $[p,q]$-coalgebras, which will form the category of morphisms from $p$ to $q$ in the underlying bicategory of $\org$.

\begin{definition}\label{coalgebras}
For polynomials $p,q\in\poly$  as in \eqref{eqn.poly_notation},
their \emph{internal hom} with respect to the tensor product $\otimes$ is the polynomial
\begin{equation}\label{eqn.internal_hom}
[p,q]\coloneqq \sum_{\phi\colon p \to q} \yon^{\sum\limits_{\;I \in p(1)} q[\phi(I)]}
\end{equation}
It can also be written $[p,q]\cong\prod_{I\in p(1)}\sum_{J\in q(1)}\prod_{j\in q[J]}\sum_{i\in p[I]}\yon$. 
\end{definition}

For intuition, a $[p,q]$-coalgebra (denoted $p \slashar q$) is a machine that outputs maps $\phi\colon p\to q$ and that inputs what \emph{flows} between them: pairs $(I,j)$ where $I\in p(1)$ is a position of $p$, which ``flows'' to $q$ as $J\coloneqq\phi(I)\in q(1)$, and $j\in q[J]$ is a direction of $q$, which ``flows'' backward to $p$ as $\phi(I,j)\in p[I]$. More precisely, using \cite[Definition 2.10]{spivak2021learners}, we define $[p,q]$-coalgebras as follows. 

\begin{definition}
The category $[p,q]\coalg$ has as objects pairs $\S = (S,\beta)$ where $S$ is a set and $\beta\colon S \to [p,q](S)$ is a function, and where a morphism from $\S$ to $\S'$ is a function $f\colon S \to S'$ making \eqref{eqn.coalg_map} commute. 
\begin{equation}\label{eqn.coalg_map}
\begin{tikzcd}
S \rar{\beta} \dar[swap]{f} & {[p,q]}(S) \dar{{[p,q]}(f)} \\
S' \rar[swap]{\beta'} & {[p,q]}(S')
\end{tikzcd}
\end{equation}
We refer to $S$ as the \emph{state set} and to each element $s\in S$ as a \emph{state}.
\end{definition}

Unwinding this definition, it is useful to break $\beta$ into two functions $\beta\coloneqq(\act^\beta,\upd^\beta)$, an \emph{action} function
\[\act^\beta\colon S\to\poly(p,q)=[p,q](1)\]
and, for each state $s \in S$, an \emph{update} function 
\[\upd^\beta_s\colon \sum_{I \in p(1)} q\left[\act^\beta_s(I)\right] \to S.\]
For a state $s\in S$ and position $I\in p(1)$ we often write $\act^\beta_s\colon p\to q$ and $\upd^\beta_s(I)\colon q[\act^\beta_s(I)]\to S.$ We may suppress the $\beta$ when it is clear from context, writing $\act_s$ and $\upd_s$. A coalgebra map $\S\to\S'$ is a function $S\to S'$ between the state sets that preserves actions and updates. 

When, for each $s \in S$, the update $\upd_s$ is the constant function sending everything to $s$, we say the coalgebra $\S$ is \emph{static}, as it remains constantly at $s$ regardless of the inputs $I \in p(1)$ and $j\in q[\act_s(I)]$ flowing between $p$ and $q$.

\begin{example}\label{ex.single_state}
A special case of a static $[p,q]$-coalgebra is given by a map $\phi \in \poly(p,q)$. For each such $\phi$, there is a coalgebra $\{\phi\}$ with a singleton state set and with $\act^\beta$ sending the point to $\phi$; we call it a \emph{singleton} coalgebra. 

A coalgebra is static iff it is the coproduct of singleton coalgebras.
\end{example}

More examples and intuition for $[p,q]$-coalgebras can be found in \cite{spivak2021learners}.

\subsection{Composition of hom-coalgebras}\label{sec.compose_hom_coalg}

We now describe how $[p,q]$-coalgebras behave like morphisms from $p$ to $q$.

\begin{proposition}\label{prop.def_org}
The categories $[p,q]\coalg$ form the hom-categories in a bicategory $\org$, which has polynomials as objects.
\end{proposition}

We use $\org$ to denote both the bicategory from \cref{prop.def_org} and the categorical operad in \cite[Definition 2.19]{spivak2021learners}, as both are derived from the monoidal double category $\org$ described in the following sections. For now, we merely present the identities and composites in this bicategory. Identities are easy: the identity object in $\org(p,p)=[p,p]\coalg$ is given by the one-state coalgebra $\idcoalg{p}$. 

The composition functor $\org(p,q)\times\org(q,r)\too\org(p,r)$ is defined as the composite:
\[[p,q]\coalg \times [q,r]\coalg \to \left([p,q] \otimes [q,r]\right)\coalg \too [p,r]\coalg,\]
where the first functor is the lax monoidality of $(-)\coalg\colon\poly \to \smcat$, as described in \cite[Proposition 2.13]{spivak2021learners}, and the second is given by applying $(-)\coalg$ to the usual ``composition'' map of internal-homs $[p,q] \otimes [q,r] \to [p,r]$ in $\poly$. Using \eqref{eqn.internal_hom} we see that on positions, this map takes the form\vspace{-.1cm}
\[\left([p,q] \otimes [q,r]\right)(1) = \poly(p,q) \times \poly(q,r) \To{\then} \poly(p,r) = [p,r](1)\]
and on directions it is given for $\phi\colon p \to q$ and $\psi\colon q \to r$ by the function
\[\bigg(\sum_{I \in p(1)} q[\phi(I)]\bigg) \times \bigg(\sum_{J \in q(1)} r[\psi(J)]\bigg) \from \sum_{I \in p(1)} r[\psi(\phi(I))]\]
which sends $(I,k)$ to $\big((I,\psi(\phi(I),k)),(\phi(I),j)\big)$. 

Concretely, the composite of a $[p,q]$-coalgebra $\S$ and a $[q,r]$-coalgebra $\S'$ is a $[p,r]$-coalgebra which we denote $\S\then\S'$  and define as follows:
\begin{itemize}
	\item its state set is given by $S \times S'$
	\item the action of the pair $(s,s')$ is given by the composite 
\[\act^{\beta\then\beta'}_{s,s'}\coloneqq(\act^\beta_s\then\act^{\beta'}_{s'})\colon p \to q \to r\]
	\item the update function of $(s,s')$ is induced by the functions
\begin{align*}
	\sum_{I \in p(1)} r\left[\act^{\beta\then\beta'}_{s,s'}(I)\right] \To{(I,k)\mapsto\left(I,\act^{\beta'}_{s'}\left(\act^\beta_s(I),k\right)\right)} \sum_{I \in p(1)} q\left[\act^\beta_s(I)\right] \To{\upd^\beta_s} S,\\
	\sum_{I \in p(1)} r\left[\act^{\beta\then\beta'}_{s,s'}(I)\right] \To{(I,k)\mapsto\left(\act^\beta_s(I),k\right)} \sum_{J \in q(1)} r\left[\act^{\beta'}_{s'}(J)\right] \To{\upd^{\beta'}_{s'}} S'.
\end{align*}
\end{itemize}

Horizontal composition of coalgebra-morphisms---i.e.\ of the 2-cells in the bicategory $\org$---is given simply by the cartesian product. The coherence isomorphisms and axioms for a bicategory then follow from the essential uniqueness of finite products of sets, and the unitality and associativity of composition for polynomial maps.

\subsection{Monoidal product of coalgebras}

It is shown in \cite[Proposition 2.13]{spivak2021learners} that the tensor product $\otimes$ of polynomials extends to make $\org$ a monoidal bicategory. That is, for polynomials $p,q,p',q'$ there is a functor
\[[p,q]\coalg \times [p',q']\coalg \to \left([p,q] \otimes [p',q']\right)\coalg \to [p \otimes p',q \otimes q']\coalg\]
derived from the map of polynomials $[p,q] \otimes [p',q'] \to [p {\otimes} p',q {\otimes} q']$ given on positions by 
\[\poly(p,q) \times \poly(p',q') \To{\otimes} \poly(p \otimes p',q \otimes q')\]
and on directions by, for $\phi\colon p \to q$ and $\phi'\colon p' \to q'$,
\[\bigg(\sum_{I \in p(1)} q[\phi_1(I)]\bigg) \times \bigg(\sum_{I' \in p'(1)} q'[\phi'_1(I')]\bigg) \fromm \sum_{(I,I') \in p(1) \times p'(1)} q[\phi_1(I)] \times q'[\phi'_1(I')]\]
sending $(I,I',j,j')$ to $(I,j,I',j')$.

Concretely, this tensor product takes a $[p,q]$-coalgebra $\S$ and a $[p',q']$-coalgebra $\S'$ to the $[p \otimes p',q \otimes q']$-coalgebra with states $S \times S'$, action
\[S \times S' \to \poly(p,q) \times \poly(p',q') \to \poly(p \otimes p',q \otimes q'),\]
and update described similarly componentwise. The tensor product of coalgebra morphisms is also given by the cartesian product of functions, and it is (very) tedious but ultimately straightforward to check that the essential uniqueness of products guarantees that $\otimes$ gives a monoidal structure on $\org$.

\subsection{$\org$ as a double category}

Defining $\org$ as a monoidal bicategory is sufficient for most of the constructions of $\org$-enriched structures in \cref{chap.org_enrich}. However, using a double category structure casting singleton coalgebras $\S_\phi\in[p,q]\coalg$ (see \cref{ex.single_state}) as morphisms $\phi\colon p\to q$ in $\poly$ facilitates our eventual definition of maps between dynamic structures. 

Specifically, the definition of $\org$ as a monoidal bicategory extends to a monoidal pseudo-double category with coalgebras as horizontal morphisms, maps in $\poly$ as vertical morphisms, and squares as in \eqref{eqn.square} given by maps of coalgebras from $\S\then\{\psi\}$ to $\{\phi\}\then\S'$. 
\begin{equation}\label{eqn.square}
\begin{tikzcd}
p \rar[slash, ""{name=S, below}]{\S} \dar[swap]{\phi} & q \dar{\psi} \\
p' \rar[slash, ""{name=T, above},swap]{\S'} & q'
\arrow[Rightarrow,shorten=5,from=S,to=T]
\end{tikzcd}
\end{equation}
The symbol $\slashar$ is intended to indicate that the map is ``dynamic'', changing in response to what flows between $p$ and $q$.

As $\{\phi\}$ and $\{\psi\}$ have only one state, and composition of coalgebras acts as the cartesian product on states, such a square amounts to a function $S \to S'$ making \eqref{eqn.coalg_square} commute:
\begin{equation}\label{eqn.coalg_square}
\begin{tikzcd}
S \rar{\beta} \dar[swap]{f} & {[p,q]}(S) \rar{\psi_\ast} & {[p,q']}(S) \ar[d, "{[p,q'](f)}" pos=.4] \\
S' \rar[swap]{\beta'} & {[p',q']}(S') \rar[swap]{\phi^\ast} & {[p,q']}(S')
\end{tikzcd}
\end{equation}

Identities and composites for these squares are determined by the bicategory structure, as this double category is a restriction in the vertical direction of the double category of lax-commuting squares in a bicategory.\footnote{It should be noted however that the vertical arrows in $\org$ are regarded as polynomial maps rather than coalgebras, so that they compose strictly unitally and associatively.}

We now proceed to discuss various categorical structures enriched in $\org$, which describe dynamical systems equipped with algebraic structure that lets us remove abstraction barriers when considering nested layers and complex arrangements of the components of the system.

\section{Dynamic structures via $\org$-Enrichment}\label{chap.org_enrich}

A monoidal double category is a viable setting for enriching various categorical structures. Intuitively, enrichment in $\org$ replaces the usual set of arrows between two objects in a categorical structure with a $[p,q]$-coalgebra for some choice of polynomials $p,q$. Therefore not only can each arrow be realized as a map of polynomials $p \to q$, but this map carries dynamics that encode how a position in $p$ and a direction in $q$ determine a transition from one arrow to another. The morphism ``reacts'' to what's flowing between $p$ and $q$.

Different situations call for different categorical structures to model their dynamics: some systems primarily involve many-to-one arrangements such as the wiring diagrams in \cref{fig.nesting}, others such as gradient descent fit naturally into a many-to-many arrow framework, and we expect in future work to consider evolving systems in which different components operate at differing time scales. Rather than choose one such categorical form to favor, and then go through the tedious exercise of forcing all of the others to conform to it, we describe how to add dynamics to the definitions of many different structures.

\slogan{A dynamic *thing* is a *thing* enriched in $\org$.}

This slogan is intentionally imprecise, so as to be maximally inclusive of different notions of categorical structures  (*things*) and notions of enrichment, and also to allow the reader who has an intuitive understanding and no need for precision to skip the remainder of this paragraph. Our intuition and examples come from the theories of enrichment described in \cite{leinster1999generalized} and \cite{shapiro2022enrichment}. In the former, a *thing* can be any suitable type of generalized multicategory, while in the latter a *thing* can be any structure defined as an algebra for a familial monad on a presheaf category equipped with a choice of ``higher'' and ``lower'' dimensional cell shapes. In both cases, *things* are algebras for a particular cartesian monad $T$ and admit an ``enriched'' analogue with respect to any $T$-multicategory. To define $T$-algebras enriched in $\org$ is then to identify $\org$ with a $T$-multicategory, and in all of our examples this identification arises naturally from the observation that monoidal double categories give rise to $T$-multicategories in a natural way.

We now give specific instances of $\org$-enrichment: in \cref{sec.org_cats} for dynamical categories, in \cref{sec.org_multicats} for dynamical multicategories and operads, and in \cref{sec:org_monoidalcats} for dynamical monoidal categories and PRO(P)s. We are also interested in using dynamic duoidal categories to describe dynamical systems in which different contributors to a system operate at different rates, using the duoidal structure on $\org$ based on $\tri$, but that is beyond the scope of this paper.

\subsection{Dynamic categories}\label{sec.org_cats}

Enrichment of categories only uses $\org$'s double category structure---not its monoidal structure---as any double category forms an $f\!c$-multicategory (also known as a virtual double category) in the sense of \cite[Section 2.1]{leinster1999generalized}. 
The following definition of enrichment in $\org$ is an unwound version of the more general definition in \cite[Section 2.2]{leinster1999generalized}.

\begin{definition}\label{def.org_enriched_cat}
An $\org$-enriched (henceforth \emph{dynamic}) category $A$ consists of
\begin{itemize}
	\item a set $A_0$ of objects;
	\item for each $a \in A_0$, a polynomial $p_a$;
	\item for each $a,b \in A_0$, a $[p_a,p_b]$-coalgebra $\S_{a,b}$;
	\item for each $a \in A_0$, an ``identitor'' square in $\org$ as in \eqref{eqn.tors} left; and
	\item for each $a,b,c \in A_0$, a ``compositor'' square in $\org$ as in \eqref{eqn.tors} right:
\begin{equation}\label{eqn.tors}
\begin{tikzcd}
p_a \dar[equals] \rar[slash, ""{name=S, below}]{\idcoalg{p_a}} & p_a \dar[equals] \\
p_a \rar[slash, ""{name=T, above},swap]{\S_{a,a}} & p_a
\arrow[Rightarrow,shorten=5,from=S,to=T]
\end{tikzcd}\qquad\qquad\begin{tikzcd}
p_a \dar[equals] \rar[slash]{\S_{a,b}} & p_b \rar[slash]{\S_{b,c}} & p_c \dar[equals] \\
p_a \ar[slash, ""{name=T, above}]{rr}[swap]{\S_{a,c}} & & p_c
\arrow[Rightarrow,shorten=4,from=1-2,to=T]
\end{tikzcd}
\end{equation}
\end{itemize}
such that these squares satisfy unit and associativity equations (\cref{CatEquations}).
\end{definition}

The sets $S_{a,b}$ form an ordinary category which we say \emph{underlies} $A$. 
In fact, a dynamic category could be equivalently defined as an ordinary category such that each object $a$ is assigned a polynomial $p_a$ and each set of arrows $\Hom(a,b)$ is equipped with a $[p_a,p_b]$-coalgebra structure, with composition and identities respecting the coalgebra structure. This means that the arrow $\id_a$ in $\Hom(a,a)$ acts as the identity map on $p_a$ and is unchanged by updates, while for $f$ in $\Hom(a,b)$ and $g$ in $\Hom(b,c)$ the composite $f\then g$ acts as the composite $p_a \to p_b \to p_c$ of the actions of $f$ and $g$, and the update of their composite equals the composite of their updates.

\subsection{Dynamic multicategories}\label{sec.org_multicats}

A monoidal double category also gives rise to an $f\!m$-multicategory in the sense of \cite[Section 3.1]{leinster1999generalized}, 
so we can talk about multicategories enriched in $\org$ as in \cite[Section 3.2]{leinster1999generalized}.

\begin{definition}
An $\org$-enriched (henceforth \emph{dynamic}) multicategory $A$ consists of
\begin{itemize}
	\item a set $A_0$ of objects;
	\item for each $a \in A_0$, a polynomial $p_a$;
	\item for each $a_1,...,a_n,b \in A_0$, a $[p_{a_1} \otimes \cdots \otimes p_{a_n},p_b]$-coalgebra $\S_{a_1,...,a_n\,;\,b}$;
	\item for each $a \in A_0$, an ``identitor'' square in $\org$ as in \eqref{eqn.multi_tors} left; and
	\item for each $a_{1,1},\ldots,a_{1,m_1},\;\ldots,\;a_{n,1},\ldots,a_{n,m_n},\;b_1,\ldots,b_n,$ and $c \in A_0$, a ``compositor'' square in $\org$ as in \eqref{eqn.multi_tors} right
\end{itemize}
\begin{equation}\label{eqn.multi_tors}
  \begin{tikzcd}[ampersand replacement=\&]
  p_a \dar[equals] \rar[slash, ""{name=S, below}]{\idcoalg{p_a}} \& p_a \dar[equals] \\
  p_a \rar[slash, ""{name=T, above},swap]{\S_{a;a}} \& p_a
  \arrow[Rightarrow,shorten=5,from=S,to=T]
  \end{tikzcd}
 \qquad\quad
  \begin{tikzcd}[column sep=huge, ampersand replacement=\&]
  p_{a_{1,1}} \otimes \cdots \otimes p_{a_{n,m_n}} \dar[equals] \ar[r, slash, "\bigotimes_i \S_{a_{i,1},...,a_{i,m_i};b_i}", ""' name=U] \&[20pt] p_{b_1} \otimes \cdots \otimes p_{b_n} \rar[slash]{\S_{b_1,...,b_n\,;\,c}} \&[-10pt] p_c \dar[equals] \\
  p_{a_{1,1}} \otimes \cdots \otimes p_{a_{n,m_n}} \ar[slash, ""{name=T}]{rr}[swap]{\S_{a_{1,1},...,a_{n,m_n};c}} \& \& p_c
  \arrow[Rightarrow,shorten=6,from=U-|T,to=T]
  \end{tikzcd}
\end{equation}
satisfying unit and associativity equations (see \cref{operadequations} for the one-object case).
\end{definition}

The sets $S_{a,b}$ form an ordinary (set-enriched) multicategory, which underlies $A$ and has a description similar to the underlying category we described below \cref{def.org_enriched_cat}. 

We will mostly be interested in what we call a \emph{dynamic operad}, the case when a dynamic multicategory $A$ has only one object, assigned the polynomial ``interface'' $p$. It consists simply of a $[p^{\otimes n},p]$-coalgebra $\S_n$ for each $n \in \nn$, equipped with coalgebra maps
\begin{equation}\label{eqn.org_operad}
\idcoalg{p} \to \S_1
\qqand
\bigotimes_{i\in I} \S_{n_i} \to \S_N
\end{equation}
where $I$ is any finite set and $N\coloneqq\sum_{i\in I}n_i$, which together satisfy the usual equations. 

\begin{example}\label{ex.collective}
A \emph{collective} (as defined in \cite{niu2021collectives}) is a $\otimes$-monoid in $\poly$, meaning a polynomial $p$ equipped with a monoid structure on its positions $p(1)$ and co-unital co-associative ``distribution'' functions $p[I \cdot J] \to p[I] \times p[J]$ for each $I,J \in p(1)$. This can be viewed as a dynamic operad where $\S_n$ is given by $\{\cdot_n\}$, the singleton coalgebra on the $n$-ary monoidal product $(\cdot_n)\colon p^{\otimes n} \to p$, and where the maps of coalgebras in \eqref{eqn.org_operad} are isomorphisms deduced from the equations for a monoid.
\end{example}

\begin{example}
In \cref{ex.collective}, the coalgebras $\S_n$ are determined by a single map of polynomials, with trivial updates since the state sets are singletons. This can be generalized to an intermediate notion between collectives and dynamic multicategories, where the coalgebras are still static but may have multiple states.

Given any multicategory $M$ and multifunctor $F\colon M \to\poly$, where $\poly$ here denotes the multicategory underlying $(\poly,\yon,\otimes)$, there is a dynamic multicategory $A_F$ with 
\begin{itemize}
	\item object set $\ob(M)$;
	\item for each $a \in \ob(M)$, the polynomial interface $p_a \coloneqq F(a)$;
	\item for each tuple $(a_1,...,a_n\,;\,b)$ in $\ob(M)$, state set $S_{a_1,...,a_n\,;\,b} \coloneqq M(a_1,...,a_n\,;\,b)$;
	\item the action $\act^\beta\colon M(a_1,...,a_n\,;\,b) \to \poly(p_{a_1} \otimes \cdots \otimes p_{a_n},p_b)$ is given by $F$; and 
	\item for any state $s$ in $M(a_1,...,a_n\,;\,b)$, the update function $\upd^\beta_s$ is the constant function at $s$.
	\qedhere
\end{itemize}
\end{example}

\begin{example}
Let $\S$ be any $p$-coalgebra for a polynomial $p$. There is a dynamic operad on $p$ with $\S_0\coloneqq \S$, with $\S_1\coloneqq\idcoalg{p}$, and with all other $\S_n\coloneqq\varnothing$ assigned the empty coalgebra.
\end{example}

\begin{example}
Consider a dynamic operad with interface $\yon\in\poly$. The internal hom polynomial $[\yon^{\otimes n},\yon]$ is simply $\yon$, so this structure amounts to an operad $S$ with a function $S_n \to S_n$ for each $n$, commuting with the operad structure. A dynamic operad on $\yon$ can thus be identified with an operad $\cat{S}$ equipped with an operad map $\cat{S}\to\cat{S}$ to itself.
\end{example}

\subsection{Dynamic monoidal categories}\label{sec:org_monoidalcats}

A monoidal double category is precisely a representable $f\!m\!c$-multicategory as in \cite[Section 2]{shapiro2022enrichment}, 
so we can also enrich strict monoidal categories in $\org$.\footnote{We use throughout the notion of \emph{strong} enrichment in a monoidal double category from \cite[Section 3]{shapiro2022enrichment}.} These are similar to $\org$-enriched multicategories, but include many-to-many coalgebras rather than just many-to-one. 

\begin{definition}\label{enriched_monoidal}
An $\org$-enriched (henceforth \emph{dynamic}) strict monoidal category $A$ consists of
\begin{itemize}
	\item a monoid $(A_0,e,*)$ of objects;
	\item for each $a \in A_0$, a polynomial $p_a$;
	\item an isomorphism of polynomials $y \cong p_e$;
	\item for each $a,a' \in A_0$, an isomorphism of polynomials $p_{a} \otimes p_{a'} \cong p_{a*a'}$;
	\item for each $a,b \in A_0$, a $[p_a,p_b]$-coalgebra $\S_{a,b}$;
	\item for each $a \in A_0$, an ``identitor'' square in $\org$ as in \cref{eqn.adaptive_tor} left;
	\item for each $a,b,c \in A_0$, a ``compositor'' square in $\org$ as in \cref{eqn.adaptive_tor} center; and
	\item for each $a,a',b,b' \in A_0$, a ``productor'' square in $\org$ as in \cref{eqn.adaptive_tor} right:
\end{itemize}
\begin{equation}\label{eqn.adaptive_tor}
\begin{tikzcd}[column sep=35pt]
p_a \dar[equals] \rar[slash, ""{name=S, below}]{\idcoalg{p_a}} & p_a \dar[equals] \\
p_a \rar[slash, ""{name=T, above},swap]{\S_{a,a}} & p_a
\arrow[Rightarrow,shorten=5,from=S,to=T]
\end{tikzcd}
\qquad
\begin{tikzcd}[column sep=30pt]
p_a \dar[equals] \rar[slash]{\S_{a,b}} & p_b \rar[slash]{\S_{b,c}} & p_c \dar[equals] \\
p_a \ar[slash, ""{name=T, above}]{rr}[swap]{\S_{a,c}} & & p_c
\arrow[Rightarrow,shorten=4,from=1-2,to=T]
\end{tikzcd}
\qquad
\begin{tikzcd}[column sep=50pt]
p_a \otimes p_{a'} \dar[equals,swap]{\wr} \rar[slash,""{name=S,below}]{\S_{a,b} \otimes \S_{a',b'}} & p_b \otimes p_{b'} \dar[equals,swap]{\wr} \\
p_{a*a'} \rar[slash, ""{name=T, above},swap]{\S_{a*a',b*b'}} & p_{b*b'}
\arrow[Rightarrow,shorten=5,from=S,to=T]
\end{tikzcd}
\end{equation}
satisfying unit, associativity, and interchange equations (see \cref{PROequations} for the one-object case).
\end{definition}

Similar to \cref{sec.org_cats,sec.org_multicats}, the sets $S_{a,b}$ form the arrows in an ordinary strict monoidal category underlying $A$.  

For the rest of this paper, we will only be interested in the restricted case of a dynamic monoidal category with object monoid $(\nn,0,+)$, which we call a dynamic PRO.%
\footnote{A PRO is the non-symmetric version of a PROP. While all of our examples are in fact symmetric, to keep the paper short we do not describe their symmetry operations.} 
Concretely, this consists of a polynomial interface $p$ (so that in the notation above $p_n\coloneqq p^{\otimes n}$ for $n \in \nn$) along with a $[p^{\otimes m},p^{\otimes n}]$-coalgebra $S_{m,n}$ for each $m,n \in \nn$, equipped with the maps of coalgebras as in \eqref{eqn.adaptive_tor}. The identitors, compositors, productors, and their equations amount to the ability to compose any string diagram of the usual sort for monoidal categories, with each $m$-to-$n$ box given by a state in $S_{m,n}$, into a new box (i.e.\ state) with the appropriate sources and targets. We denote a dynamic PRO as $(p,\S)$, where $\S$ encodes all of the coalgebras $\S_{m,n}$ that constitute the $\org$-enrichment and the structure maps are implicit.

\bigskip
We now turn to morphisms between dynamic PROs; the interested reader can hopefully find analogous definitions for dynamic categories and operads.

\begin{definition}
A \emph{morphism} of dynamic PROs from $(p,\S)$ to $(p',\S')$ is given by a map of polynomials $\phi\colon p \to p'$ and, for each $m,n \in \nn$, ``commutor'' squares as in \eqref{eqn.adaptive_map} in $\org$ which commute with the identitor, compositor, and productor squares.
\setlength{\belowdisplayskip}{-5pt}
\begin{equation}\label{eqn.adaptive_map}
\begin{tikzcd}
p^{\otimes m} \rar[slash, ""{name=S, below}]{\S_{m,n}} \dar[swap]{\phi^{\otimes m}} & p^{\otimes n} \dar{\phi^{\otimes n}} \\
p'^{\otimes m} \rar[slash, ""{name=T, above}, swap]{\S'_{m,n}} & p'^{\otimes n}
\arrow[Rightarrow,shorten=5,from=S,to=T]
\end{tikzcd}
\end{equation}
\end{definition}
\setlength{\belowdisplayskip}{11pt}

This definition of morphism (taken from \cite[Section 3]{shapiro2022enrichment}) 
is the direct theoretical benefit of treating $\org$ as a monoidal double category rather than as a monoidal bicategory (closer to its description in \cite{spivak2021learners}). Otherwise morphisms could either only be easily defined between dynamic PROs with the same interface polynomial, which is needlessly restrictive, or take the form of a $[p,p']$-coalgebra, which seems to be too general to be of much use.

A morphism $(p,\S) \to (p',\S')$ can be interpreted as a way of telling the codomain how to run the domain. The map of polynomials $p \to p'$ specifies how the positions of $p$ can be modeled by those of $p'$ and how the directions of $p'$ are returned as directions of $p$, while the commutor squares describe how the states of $\S_{m,n}$ can be modeled by those of $\S'_{m,n}$ in a way that respects this change of interface. A type of theorem that we hope to instantiate in future work is of the form ``this dynamic structure that we're interested in can be run by (has a map to) this other dynamic structure that we already understand well.''

\begin{example}
For a fixed polynomial $p$, there is a terminal dynamic PRO with interface $p$, which we denote $\S^{p!}$; here $\S^{p!}_{m,n}$ is the terminal $[p^{\otimes m},p^{\otimes n}]$-coalgebra for each $m,n\in\nn$. 

A state in $\S^{p!}$ is a (not necessarily finite) $[p^{\otimes m},p^{\otimes n}]$-tree. By this we mean a tree co-inductively defined by a root node labeled with a polynomial map $\phi\colon p^{\otimes m} \to p^{\otimes n}$ together with an arrow---whose source is the root and whose target is another $[p^{\otimes m},p^{\otimes n}]$-tree---assigned to each tuple 
\begin{equation}\label{eqn.pmn_directions}
\big((I_1,\ldots,I_m), i_1,\ldots,i_n\big) \in  p^{\otimes m}(1)\times p^{\otimes n}[\phi(I_1,...,I_m)] 
\end{equation}

The action of such a tree is simply the map $\phi$ labeling its root, and the update sends a tuple as in \eqref{eqn.pmn_directions} to the target of its assigned arrow. 

The idea is that the state-set of the terminal dynamic PRO encodes all possible trajectories along different actions, and this coalgebra is terminal because from any other coalgebra there is a map to $\S^{p!}_{m,n}$ sending each state to the tree whose root is labeled by the action of the state and whose edges from the root go to the trees for each of the state's possible updates.

To define a dynamic PRO structure on the terminal coalgebra $\S^{p!}$, it only remains to define maps of coalgebras as in \cref{eqn.adaptive_tor}, and these are all taken to be the unique map to the terminal $[p^{\otimes m},p^{\otimes n}]$-coalgebra; the equations hold automatically. This is the terminal dynamic PRO with interface $p$ because for any other such dynamic PRO there is a morphism given by the identity map on $p$ and with commutor squares to $\S^{p!}_{m,n}$ the unique map to the terminal $[p^{\otimes m},p^{\otimes n}]$-coalgebra. In other words, $\S^{p!}$ \emph{uniquely runs} any other dynamic PRO with interface $p$.
\end{example}

\section{Dynamic Structures in Nature}

Our main results are that dynamic structures describe phenomena we see instantiated around us. In this paper, we focus on deep learning and a prediction market in which the reputations of various guess-makers evolve based on how successful they are.

\subsection{The prediction market dynamic operad}\label{sec.kelley}

Fix a finite set $X$, elements of which we call \emph{outcomes} and intuit to be ``all equally likely'', define the set $\bet_X$ of \emph{guesses on $X$} as\footnote{We assume that each guess assigns a nonzero probability to each possible outcome, which avoids the issues of dividing by zero when updating or permanent loss of a guess-maker's reputation. This should be interpreted as both humility and good strategy on the part of the guess-makers.}
\[
	\bet_X\coloneqq\left\{\gamma\colon X\to(0,1]\;\;\middle|\;\;1=\sum_x\gamma(x) \right\}
\]
Let $\Delta^+$ denote the operad of finite nowhere-zero probability distributions, where $\Delta^+_N$ is defined as above with the natural number $N$ regarded as the $N$-element set. 
Then $\bet_X$ is an algebra for it: for any $\mu\in\Delta_N$ and $\gamma \in (\bet_X)^N$, we define 
\[
	\mu\cdot\gamma\coloneqq\bigg(x\mapsto\sum_{i\in N}\mu_i\cdot\gamma_i(x)\bigg)
\]
and it is easy to check that $(\mu\cdot\gamma)\in\bet_X$, i.e.\ its components are in bounds $(\mu\cdot\gamma)(x)\in (0,1]$ and it is normalized $\sum_x(\mu\cdot\gamma)(x)=1$.

We now construct a dynamic operad with interface $p_X\in\poly$ defined as:
\[
p_X\coloneqq \bet_X\,\yon^X
\]
and use the $\Delta^+_N$ as our state spaces. The idea is that a state $\mu\in\Delta^+_N$ says how much the organization trusts each of its $N$ members (guess-makers) relative to each other. A member's position at a given moment is a report of how much confidence it has in each of the $X$-many possibilities, represented by its probability distribution.

The action of a trust distribution $\mu \in \Delta^+_N$ is the map of polynomials $p_X^{\otimes N} \to p_X$ which on positions sends $\gamma \in (\bet_X)^N$ to $\mu \cdot \gamma$ and on directions sends $x \in X$ to $(x,...,x) \in X^N$. The idea is that the organization aggregates its members' predictions according to its current trust-distribution, and the outcome is accurately communicated back to each member.

The most interesting part of the dynamic structure is how the trust distribution is updated once predictions are made and a result $x\in X$ is returned. When $N=0$, there's nothing to do: $\Delta^+_0=\varnothing$. For membership $N\geq 1$, trust distribution $\mu\in\Delta^+_N$, guesses $\gamma\in(\bet_X)^N$, and outcome $x\in X$, we define the updated trust distribution $\gamma(x) * \mu \in\Delta^+_N$ as
\[
\gamma(x) * \mu \coloneqq \left( i \mapsto \frac{\gamma_i(x)\mu_i}{\sum_j \gamma_j(x)\mu_j}\right).
\]

Finally, we describe the operadic structure maps. As $\Delta^+_1$ is a singleton set whose action is the identity on $p_X$, the identitor $\{\id_{p_X}\} \to \Delta^+_1$ is an isomorphism. The operadic compositor is given by the usual operad structure on (nowhere-zero) distributions:
\[
\Delta^+_N \times \Delta^+_{M_1} \times \cdots \times \Delta^+_{M_N} \to \Delta^+_{\sum_i M_i} \qquad\qquad (\mu,\nu_1,\ldots,\nu_N) \mapsto \mu \circ \nu \coloneqq \left( (i,j) \mapsto \mu_i\nu_j \right).
\]

\begin{theorem}\label{predictionadaptive}
The maps defined above are maps of coalgebras and satisfy the coherence equations of a dynamic operad described in \cref{operadequations}.
\end{theorem}

This is proven in \cref{proofs}.

\subsection{The gradient descent dynamic PRO}

Deep learning uses the algorithm of gradient descent to optimize a choice of function, based on external feedback on its output. This naturally fits into the paradigm of dynamic structures, since functions $\rr^m \to \rr^n$ can form the states of a polynomial coalgebra, with the feedback providing the information needed to update the choice of function. These functions can be composed and juxtaposed in a way that preserves the updates. That is, the composite of gradient descenders is a gradient descender.

\begin{definition}\label{def.Smn}
For the rest of this section, we will use the state sets 
\[
S_{m,n} \coloneqq \left\{(M \in \nn, f\colon \rr^{M+m} \to \rr^n, p \in \rr^M) \;\middle|\; f \textrm{ is differentiable}\right\}.
\qedhere
\]
\end{definition}

The idea is that these states are the possible parameters among which a gradient descender is meant to find the optimal choice, while $f$ dictates how the parameter affects the resulting function $f(p,-)$. In the dynamics of these states described below, only the value of the parameter $p$ will be updated; the parameter-space dimension $M$ and the parameterized function $f$ will remain fixed, though network composition of gradient descenders will involve combining these data in nontrivial ways. Fix $\epsilon>0$.

For every $x\in\rr$, let $T_x\rr$ denote the tangent space at $x$; for all practical purposes $T_x \rr$ can be regarded as simply $\rr$, but in both the description of polynomials as bundles and the intuition for this example it makes sense to use the tangent space. We proceed to define a dynamic PRO with interface $t \coloneqq \sum_{x \in \rr} \yon^{T_x \rr}$ and coalgebras $\S_{m,n}$ which update the state sets $S_{m,n}$ from \cref{def.Smn} using gradient descent. The PRO structure maps encode how networks of gradient descenders can be composed into a single gradient descender with a larger parameter space.

\begin{definition}
The $[t^{\otimes m},t^{\otimes n}]$-coalgebra structure on $S_{m,n}$ is given by 
\begin{itemize}
	\item On positions, the action $\act^\beta_{M,f,p}\colon \rr^m \to \rr^n$ is given by $f(p,-)$.
	\item For $x \in \rr^m$, the action $\act^\beta_{M,f,p}(x,-)\colon T_{f(p,x)} \rr^n \to T_x \rr^m$ on directions sends $y\in T_{f(p,x)}$ to $\pi_m (Df)^\top \cdot y$.
	\item The update function $\upd^\beta_{M,f,p}$ sends $x \in \rr^m$ and $y \in T_{f(p,x)}$ to $(M,f,p+\epsilon \pi_M (Df)^\top \cdot y)$ for our fixed $\epsilon$.
	\qedhere
\end{itemize}
\end{definition}

The action of a state as a map $t^{\otimes m} \to t^{\otimes n}$ is given by applying the parameterized function $f$ with the parameter $p$, resulting in a function $\rr^m \to \rr^n$ as desired. The transpose $(Df)^\top$ of the derivative of $f$ sends a feedback vector $y \in T_{f(p,x)} \rr^n$, which can be interpreted as the difference in $\rr^n$ between the ``correct'' result for $x$ and the current approximation $f(p,x)$, to the corresponding ``correction'' to $(p,x)$ in $\rr^{M+m}$. The projection of this correction to $T_x \rr^m$ provides the action of the state on directions, which in a network will then be further propagated back to the gradient descender which had output $x$. The projection to $T_p \rr^M$ provides the direction and magnitude in which to update the parameters (scaled by the ``learning rate'' $\epsilon$).

Thus far, we have provided the data of the polynomial $t$ and the $[t^{\otimes m},t^{\otimes n}]$-coalgebras $\S_{m,n}$ needed to define a dynamic PRO. We now define the identitor, compositor, and productor morphisms of coalgebras presented by the squares in \cref{enriched_monoidal}.
\begin{itemize}
	\item The identitors $\idcoalg{t^{\otimes n}} \to \S_{n,n}$ send the unique state in the domain to the state 
\[(0,\id_{\rr^n},0) \in S_{n,n}.\] 
	\item The compositors $\S_{\ell,m}\then\S_{m,n} \to \S_{\ell,n}$ send the pair $((L,f,p),(M,g,q))$ to 
\[\left( M+L,\,g(-,f(-,-))\colon \rr^{M+L+\ell} \To{\id \times f} \rr^{M+m} \To{g} \rr^n,\, (q,p) \in \rr^{M+L} \right).\]
	\item The productors $\S_{m,n} \otimes \S_{m',n'} \to \S_{m+m',n+n'}$ send the pair $((M,f,p),(M',f',p'))$ to 
\[(M+M',\,(f,f'),\,(p,p')).\]
\end{itemize}

These structure maps ensure that whenever two gradient descenders are combined in series or parallel, the resulting composite descender retains the parameter spaces of both. Likewise when the input or output of a descender is wired past some other descender in a network, it does not contribute any new parameters and merely preserves its input/output until plugged into a descender. The following is proven in \cref{proofs}.

\begin{theorem}\label{gradientadaptive}
The maps defined above are maps of coalgebras and satisfy the coherence equations of a dynamic PRO described in \cref{PROequations}.
\end{theorem}

\appendix
\section{Coherence Equations}\label{coherences}

We now present the equations that must be satisfied by the structure maps in dynamic categories, operads and PROs. While we only provide the equations for the single-object variant of dynamic multicategories and monoidal categories, respectively, the equations in the general case are entirely analogous.

\begin{definition}\label{CatEquations}
The equations between the identitors and compositors in a dynamic category are as follows:
\begin{itemize}
	\item The left and right unit laws
\begin{equation}\label{eqn.cat_unit}
\begin{tikzcd}[ampersand replacement=\&]
  p_a \dar[equals] \ar[r, slash, "\{\id_{p_a}\}", ""' name=R] \& p_a \rar[slash, ""' name=S]{\S_{a,b}} \dar[equals] \& p_b \dar[equals] \\
  p_a \dar[equals] \ar[r,""' name=U,""{name=W,below},slash,swap,"\S_{a,a}"] \& p_a \ar[r,""' name=V,slash,swap,"\S_{a,b}"] \& p_b \dar[equals] \\
  p_a \ar[slash, ""{name=T}]{rr}[swap]{\S_{a,b}} \& \& p
  \arrow[Rightarrow,shorten=8,from=W-|T,to=T]
  \arrow[Rightarrow,shorten=7,from=R,to=U]
  \arrow[equals,shorten=8,from=S,to=V]
  \end{tikzcd} \quad = \quad \begin{tikzcd}[ampersand replacement=\&]
p_a \dar[equals] \rar[slash,""{name=S, below}]{\S_{a,b}} \& p_b \dar[equals] \\
p_a \rar[slash, ""{name=T, above},swap]{\S_{a,b}} \& p_b
\arrow[equals,shorten=7,from=S,to=T]
  \end{tikzcd} \quad = \quad \begin{tikzcd}[ampersand replacement=\&]
  p_a \dar[equals] \ar[r, slash, "\S_{a,b}", ""' name=R] \& p_b \rar[slash, ""' name=S]{\{\id_{p_b}\}} \dar[equals] \& p_b \dar[equals] \\
  p_a \dar[equals] \ar[r,""' name=U,""{name=W,below},slash,swap,"\S_{a,b}"] \& p_b \ar[r,""' name=V,slash,swap,"\S_{b,b}"] \& p_b \dar[equals] \\
  p_a \ar[slash, ""{name=T}]{rr}[swap]{\S_{a,b}} \& \& p
  \arrow[Rightarrow,shorten=8,from=W-|T,to=T]
  \arrow[equals,shorten=8,from=R,to=U]
  \arrow[Rightarrow,shorten=7,from=S,to=V]
  \end{tikzcd}
\end{equation}
	\item The associativity law
\begin{equation}\label{eqn.cat_assoc}
\begin{tikzcd}[ampersand replacement=\&]
  p_a \dar[equals] \ar[r, slash, "\S_{a,b}", ""' name=U] \& p_b \ar[r, slash, "\S_{b,c}"] \& p_c \dar[equals] \rar[slash, ""' name=V]{\S_{c,d}} \& p_d \dar[equals] \\
  p_a \dar[equals] \ar[slash, ""{name=W,above},""{name=S, below}]{rr}[swap]{\S_{a,c}} \& \& p_c \ar[slash,""' name=X]{r}[swap]{\S_{c,d}} \& p_d \dar[equals] \\
  p_a \ar[slash, ""{name=T}]{rrr}[swap]{\S_{a,d}} \& \& \& p_d
  \arrow[Rightarrow,shorten=8,from=U-|W,to=W]
  \arrow[equals,shorten=8,from=V,to=X]
  \arrow[Rightarrow,shorten=8,from=S-|T,to=T]
  \end{tikzcd} \quad = \quad \begin{tikzcd}[ampersand replacement=\&]
  p_a \dar[equals] \ar[r, slash, "\S_{a,b}", ""' name=U] \& p_b \dar[equals] \ar[r, slash, "\S_{b,c}"] \& p_c \rar[slash, ""' name=V]{\S_{c,d}} \& p_d \dar[equals] \\
  p_a \dar[equals] \ar[slash, ""{name=W,above}]{r}[swap]{\S_{a,b}} \& p_b \ar[slash, ""{name=X,above},""{name=S, below}]{rr}[swap]{\S_{b,d}} \& \& p_d \dar[equals] \\
  p_a \ar[slash, ""{name=T}]{rrr}[swap]{\S_{a,d}} \& \& \& p_d
  \arrow[equals,shorten=8,from=U-|W,to=W]
  \arrow[Rightarrow,shorten=8,from=V-|X,to=X]
  \arrow[Rightarrow,shorten=8,from=S-|T,to=T]
  \end{tikzcd}
\end{equation}
\qedhere
\end{itemize}
\end{definition}

We now present the equations for dynamic operads. These equations derive directly from the definition of operads, namely the associativity and unitality of operadic composition, but unlike the equations above only involve a single polynomial $p$.

\begin{definition}\label{operadequations}
The equations between the identitors and compositors in a dynamic operad are as follows:
\begin{itemize}
	\item The left and right unit laws
\begin{equation}\label{eqn.operad_unit}
\begin{tikzcd}[ampersand replacement=\&]
  p^{\otimes n} \dar[equals] \ar[r, slash, "\{\id_p\}^{\otimes n}", ""' name=R] \& p^{\otimes n} \rar[slash, ""' name=S]{\S_n} \dar[equals] \& p \dar[equals] \\
  p^{\otimes n} \dar[equals] \ar[r,""' name=U,""{name=W,below},slash,swap,"\S_1^{\otimes n}"] \& p^{\otimes n} \ar[r,""' name=V,slash,swap,"\S_n"] \& p \dar[equals] \\
  p^{\otimes n} \ar[slash, ""{name=T}]{rr}[swap]{\S_n} \& \& p
  \arrow[Rightarrow,shorten=8,from=W-|T,to=T]
  \arrow[Rightarrow,shorten=7,from=R,to=U]
  \arrow[equals,shorten=8,from=S,to=V]
  \end{tikzcd} \;\;\; = \;\;\; \begin{tikzcd}[ampersand replacement=\&]
p^{\otimes n} \dar[equals] \rar[slash,""{name=S, below}]{\S_n} \& p \dar[equals] \\
p^{\otimes n} \rar[slash, ""{name=T, above},swap]{\S_n} \& p
\arrow[equals,shorten=7,from=S,to=T]
  \end{tikzcd} \;\;\; = \;\;\; \begin{tikzcd}[ampersand replacement=\&]
  p^{\otimes n} \dar[equals] \ar[r, slash, "\S_n", ""' name=R] \& p \rar[slash, ""' name=S]{\{\id_p\}} \dar[equals] \& p \dar[equals] \\
  p^{\otimes n} \dar[equals] \ar[r,""' name=U,""{name=W,below},slash,swap,"\S_n"] \& p \ar[r,""' name=V,slash,swap,"\S_1"] \& p \dar[equals] \\
  p^{\otimes n} \ar[slash, ""{name=T}]{rr}[swap]{\S_n} \& \& p
  \arrow[Rightarrow,shorten=8,from=W-|T,to=T]
  \arrow[equals,shorten=8,from=R,to=U]
  \arrow[Rightarrow,shorten=7,from=S,to=V]
  \end{tikzcd}
\end{equation}
	\item The associativity law
\[ 
\begin{tikzcd}[column sep=large,ampersand replacement=\&]
  p^{\otimes \ell_{1,1}} \otimes \cdots \otimes p^{\otimes \ell_{n,m_n}} \dar[equals,swap,"\wr"] \ar[r, slash, "\bigotimes_{i,j} \S_{\ell_{i,j}}", ""' name=U] \& p^{\otimes m_1} \otimes \cdots \otimes p^{\otimes m_n} \ar[r, slash, "\bigotimes_i \S_{m_i}"] \& p^{\otimes n} \dar[equals] \rar[slash, ""' name=V]{\S_n} \& p \dar[equals] \\
  p^{\otimes (\sum_j \ell_{1,j})} \otimes \cdots \otimes p^{\otimes (\sum_j \ell_{n,j})} \dar[equals,swap,"\wr"] \ar[slash, ""{name=W,above},""{name=S, below}]{rr}[description]{\bigotimes_i \S_{\sum_j \ell_{1,j}}} \& \& p^{\otimes n} \rar[slash,swap,""' name=X]{\S_n} \& p \dar[equals] \\
  p^{\otimes (\sum_{i,j} \ell_{i,j})} \ar[slash, ""{name=T}]{rrr}[swap]{\S_{\sum_{i,j} \ell_{i,j}}} \& \& \& p
  \arrow[Rightarrow,shorten=8,from=U-|W,to=W]
  \arrow[equals,shorten=8,from=V,to=X]
  \arrow[Rightarrow,shorten=8,from=S-|T,to=T]
  \end{tikzcd}
\]
\begin{equation}\label{eqn.operad_assoc}
=
\end{equation}
\[
\begin{tikzcd}[column sep=large,ampersand replacement=\&]
  p^{\otimes \ell_{1,1}} \otimes \cdots \otimes p^{\otimes \ell_{n,m_n}} \dar[equals] \ar[r, slash, "\bigotimes_{i,j} \S_{\ell_{i,j}}", ""' name=U] \& p^{\otimes m_1} \otimes \cdots \otimes p^{\otimes m_n} \dar[equals,swap,"\wr"] \ar[r, slash, "\bigotimes_i \S_{m_i}"] \& p^{\otimes n} \rar[slash, ""' name=V]{\S_n} \& p \dar[equals] \\
  p^{\otimes \ell_{1,1}} \otimes \cdots \otimes p^{\otimes \ell_{n,m_n}} \dar[equals,swap,"\wr"] \ar[slash, ""{name=W,above}]{r}[swap]{\bigotimes_{i,j} \S_{\ell_{i,j}}} \& p^{\otimes (\sum_i m_i)} \ar[slash, ""{name=X,above},""{name=S, below}]{rr}[swap]{\S_{\sum_i m_i}} \& \& p \dar[equals] \\
  p^{\otimes (\sum_{i,j} \ell_{i,j})} \ar[slash, ""{name=T}]{rrr}[swap]{\S_{\sum_{i,j} \ell_{i,j}}} \& \& \& p
  \arrow[equals,shorten=8,from=U-|W,to=W]
  \arrow[Rightarrow,shorten=8,from=V-|X,to=X]
  \arrow[Rightarrow,shorten=8,from=S-|T,to=T]
  \end{tikzcd}
\]
\qedhere
\end{itemize}
\end{definition}

The equations for dynamic PROs below are similarly derived from the definition of monoidal categories, namely that composition and products of arrows are associative and unital (giving the associativity and unitality equations for compositors and productors) and products are functorial (giving the interchange equations). 

\begin{definition}\label{PROequations}
The equations between the identitors, compositors, and productors in a dynamic PRO are as follows:
\begin{itemize}
	\item The identitor interchange law
\begin{equation}\label{eqn.id_inter}
\begin{tikzcd}[column sep={130,between origins}]
p^{\otimes n} \otimes p^{\otimes n'} \dar[equals] \rar[slash, ""{name=S, below}]{\idcoalg{p^{\otimes n}} \otimes \idcoalg{p^{\otimes n'}}} & 
p^{\otimes n} \otimes p^{\otimes n'} \dar[equals] \\
p^{\otimes n} \otimes p^{\otimes n'} \dar[equals,swap]{\wr} \ar[slash, ""{name=T, above}, ""{name=U, below}]{r}[description]{\S_{n,n} \otimes \S_{n',n'}} & 
p^{\otimes n} \otimes p^{\otimes n'} \dar[equals,swap]{\wr} \\
p^{\otimes (n+n')} \rar[slash, ""{name=V, above},swap]{\S_{n+n',n+n'}} & 
p^{\otimes (n+n')}
\arrow[Rightarrow,shorten=5,from=S,to=T]
\arrow[Rightarrow,shorten=5,from=U,to=V]
\end{tikzcd}\;\; = \;\;\begin{tikzcd}[column sep={130,between origins}]
p^{\otimes n} \otimes p^{\otimes n'} \dar[equals,swap]{\wr} \rar[slash, ""{name=S, below}]{\idcoalg{p^{\otimes n}} \otimes \idcoalg{p^{\otimes n'}}} & 
p^{\otimes n} \otimes p^{\otimes n'} \dar[equals,swap]{\wr} \\
p^{\otimes (n+n')} \dar[equals] \ar[slash, ""{name=T, above}, ""{name=U, below}]{r}[description]{\idcoalg{p^{\otimes (n+n')}}} & 
p^{\otimes (n+n')} \dar[equals] \\
p^{\otimes (n+n')} \rar[slash, ""{name=V, above},swap]{\S_{n+n',n+n'}} & 
p^{\otimes (n+n')}
\arrow[equals,shorten=5,from=S,to=T,swap,"\wr"]
\arrow[Rightarrow,shorten=5,from=U,to=V]
\end{tikzcd}
\end{equation}
	\item The compositor interchange law
\[
\begin{tikzcd}[column sep={110,between origins}]
p^{\otimes \ell} \otimes p^{\otimes \ell'} \dar[equals,swap]{\wr} \rar[slash, ""{name=S, below}]{\S_{\ell,m} \otimes \S_{\ell',m'}} & 
p^{\otimes m} \otimes p^{\otimes m'} \dar[equals,swap]{\wr} \rar[slash, ""{name=U, below}]{\S_{m,n} \otimes \S_{m',n'}} & 
p^{\otimes n} \otimes p^{\otimes n'} \dar[equals,swap]{\wr} \\
p^{\otimes (\ell + \ell')} \dar[equals] \rar[slash, ""{name=T, above},swap]{\S_{\ell+\ell',m+m'}} \ar[phantom,""{name=R, below}]{rr} & 
p^{\otimes (m + m')} \rar[slash, ""{name=V, above},swap]{\S_{m+m',n+n'}} & 
p^{\otimes (n + n')} \dar[equals] \\
p^{\otimes (\ell + \ell')} \ar[slash, ""{name=W, above}]{rr}[swap]{\S_{\ell+\ell',n+n'}} & & 
p^{\otimes (n + n')}
\arrow[Rightarrow,shorten=5,from=S,to=T]
\arrow[Rightarrow,shorten=5,from=U,to=V]
\arrow[Rightarrow,shorten=5,from=R,to=W]
\end{tikzcd}
\]
\begin{equation}\label{eqn.comp_inter}
= 
\end{equation}
\[
\begin{tikzcd}[column sep={110,between origins}]
p^{\otimes \ell} \otimes p^{\otimes \ell'} \dar[equals] \rar[slash, ""{name=S, below}]{\S_{\ell,m} \otimes \S_{\ell',m'}} \ar[phantom,""{name=R, below}]{rr} & 
p^{\otimes m} \otimes p^{\otimes m'} \rar[slash, ""{name=U, below}]{\S_{m,n} \otimes \S_{m',n'}} & 
p^{\otimes n} \otimes p^{\otimes n'} \dar[equals] \\
p^{\otimes \ell} \otimes p^{\otimes \ell'} \dar[equals,swap]{\wr} \ar[slash, ""{name=S, above},""{name=T, below}]{rr}[description]{\S_{\ell,n} \otimes \S_{\ell',n'}} & &
p^{\otimes n} \otimes p^{\otimes n'} \dar[equals,swap]{\wr} \\
p^{\otimes (\ell + \ell')} \ar[slash, ""{name=U, above}]{rr}[swap]{\S_{\ell+\ell',n+n'}} & & 
p^{\otimes (n + n')}
\arrow[Rightarrow,shorten=5,from=R,to=S]
\arrow[Rightarrow,shorten=5,from=T,to=U]
\end{tikzcd}
\]
	\item The compositor associativity law
\begin{equation}\label{eqn.comp_assoc}
\begin{tikzcd}
p^{\otimes k} \dar[equals] \rar[slash]{\S_{k,\ell}} \ar[phantom,""{name=S, below}]{rr} & 
p^{\otimes \ell} \rar[slash]{\S_{\ell,m}} & 
p^{\otimes m} \dar[equals] \rar[slash, ""{name=U, below}]{\S_{m,n}} & 
p^{\otimes n} \dar[equals] \\
p^{\otimes k} \dar[equals] \ar[slash, ""{name=T, above}]{rr}[swap]{\S_{k,m}} \ar[phantom,""{name=W, below}]{rrr} & &
p^{\otimes m} \rar[slash, ""{name=V, above},swap]{\S_{m,n}} & 
p^{\otimes n} \dar[equals] \\
p^{\otimes k} \ar[slash, ""{name=X, above}]{rrr}[swap]{\S_{k,n}} & & &
p^{\otimes n}
\arrow[Rightarrow,shorten=5,from=S,to=T]
\arrow[equals,shorten=5,from=U,to=V]
\arrow[Rightarrow,shorten=5,from=W,to=X]
\end{tikzcd}\;\; = \;\;\begin{tikzcd}
p^{\otimes k} \dar[equals] \rar[slash,""{name=S, below}]{\S_{k,\ell}} & 
p^{\otimes \ell} \dar[equals] \rar[slash]{\S_{\ell,m}} \ar[phantom,""{name=U, below}]{rr} & 
p^{\otimes m} \rar[slash]{\S_{m,n}} & 
p^{\otimes n} \dar[equals] \\
p^{\otimes k} \dar[equals] \rar[slash, ""{name=T, above},swap]{\S_{k,\ell}} \ar[phantom,""{name=W, below}]{rrr} & 
p^{\otimes \ell} \ar[slash, ""{name=V, above}]{rr}[swap]{\S_{\ell,n}} & &
p^{\otimes n} \dar[equals] \\
p^{\otimes k} \ar[slash, ""{name=X, above}]{rrr}[swap]{\S_{k,n}} & & &
p^{\otimes n}
\arrow[equals,shorten=5,from=S,to=T]
\arrow[Rightarrow,shorten=5,from=U,to=V]
\arrow[Rightarrow,shorten=5,from=W,to=X]
\end{tikzcd}
\end{equation}
	\item The compositor unit laws
\begin{equation}\label{eqn.comp_unit}
\begin{tikzcd}
p^{\otimes m} \dar[equals] \rar[slash,""{name=S, below}]{\idcoalg{p^{\otimes m}}} & p^{\otimes m} \dar[equals] \rar[slash,""{name=U, below}]{\S_{m,n}} & p^{\otimes n} \dar[equals] \\
p^{\otimes m} \dar[equals] \rar[slash,""{name=T, above},swap]{\S_{m,m}} \ar[phantom,""{name=W, below}]{rr} & p^{\otimes m} \rar[slash,""{name=V, above},swap]{\S_{m,n}} & p^{\otimes n} \dar[equals] \\
p^{\otimes m} \ar[slash, ""{name=X, above}]{rr}[swap]{\S_{m,n}} & & p^{\otimes n}
\arrow[Rightarrow,shorten=5,from=S,to=T]
\arrow[equals,shorten=5,from=U,to=V]
\arrow[Rightarrow,shorten=5,from=W,to=X]
\end{tikzcd}\;\; = \;\;\begin{tikzcd}
p^{\otimes m} \dar[equals] \rar[slash,""{name=S, below}]{\S_{m,n}} & p^{\otimes n} \dar[equals] \\
p^{\otimes m} \rar[slash, ""{name=T, above},swap]{\S_{m,n}} & p^{\otimes n}
\arrow[equals,shorten=5,from=S,to=T]
\end{tikzcd}\;\; = \;\;\begin{tikzcd}
p^{\otimes m} \dar[equals] \rar[slash,""{name=S, below}]{\S_{m,n}} & p^{\otimes n} \dar[equals] \rar[slash,""{name=U, below}]{\idcoalg{p^{\otimes m}}} & p^{\otimes n} \dar[equals] \\
p^{\otimes m} \dar[equals] \rar[slash,""{name=T, above},swap]{\S_{m,n}} \ar[phantom,""{name=W, below}]{rr} & p^{\otimes n} \rar[slash,""{name=V, above},swap]{\S_{n,n}} & p^{\otimes n} \dar[equals] \\
p^{\otimes m} \ar[slash, ""{name=X, above}]{rr}[swap]{\S_{m,n}} & & p^{\otimes n}
\arrow[equals,shorten=5,from=S,to=T]
\arrow[Rightarrow,shorten=5,from=U,to=V]
\arrow[Rightarrow,shorten=5,from=W,to=X]
\end{tikzcd}
\end{equation}
	\item The productor associativity law
\[
\begin{tikzcd}[column sep={180,between origins}]
p^{\otimes m} \otimes p^{\otimes m'} \otimes p^{\otimes m''}  \dar[equals,swap]{\wr} \rar[slash,""{name=S,below}]{\S_{m,n} \otimes \S_{m',n'} \otimes \S_{m'',n''}} & p^{\otimes n} \otimes p^{\otimes n'} \otimes p^{\otimes n''} \dar[equals,swap]{\wr} \\
p^{\otimes (m + m')} \otimes p^{\otimes m''} \dar[equals,swap]{\wr} \ar[slash,""{name=T,above},""{name=U,below}]{r}[description]{\S_{m+m',n+n'} \otimes \S_{m'',n''}} & p^{\otimes (n + n')} \otimes p^{\otimes n''} \dar[equals,swap]{\wr} \\
p^{\otimes (m + m' + m'')} \rar[slash, ""{name=V, above},swap]{\S_{m+m'+m'',n+n'+n''}} & p^{\otimes (n + n' + n'')}
\arrow[Rightarrow,shorten=5,from=S,to=T]
\arrow[Rightarrow,shorten=5,from=U,to=V]
\end{tikzcd}
\]
\begin{equation}\label{eqn.prod_assoc}
= 
\end{equation}
\[
\begin{tikzcd}[column sep={180,between origins}]
p^{\otimes m} \otimes p^{\otimes m'} \otimes p^{\otimes m''}  \dar[equals,swap]{\wr} \rar[slash,""{name=S,below}]{\S_{m,n} \otimes \S_{m',n'} \otimes \S_{m'',n''}} & p^{\otimes n} \otimes p^{\otimes n'} \otimes p^{\otimes n''} \dar[equals,swap]{\wr} \\
p^{\otimes m} \otimes p^{\otimes (m+m'')} \dar[equals,swap]{\wr} \ar[slash,""{name=T,above},""{name=U,below}]{r}[description]{\S_{m,n} \otimes \S_{m'+m'',n'+n''}} & p^{\otimes n} \otimes p^{\otimes (n' + n'')} \dar[equals,swap]{\wr} \\
p^{\otimes (m + m' + m'')} \rar[slash, ""{name=V, above},swap]{\S_{m+m'+m'',n+n'+n''}} & p^{\otimes (n + n' + n'')}
\arrow[Rightarrow,shorten=5,from=S,to=T]
\arrow[Rightarrow,shorten=5,from=U,to=V]
\end{tikzcd}
\]
	\item The productor unit laws
\begin{equation}\label{eqn.prod_unit}
\hspace{-1cm}\begin{tikzcd}[column sep={120,between origins}]
p^{\otimes m} \dar[equals,swap]{\wr} \rar[slash,""{name=Q, below}]{\S_{m,n}} & p^{\otimes n} \dar[equals,swap]{\wr} \\
p^{\otimes 0} \otimes p^{\otimes m} \dar[equals] \ar[slash, ""{name=S, below}, ""{name=R, above}]{r}[description]{\idcoalg{p^{\otimes 0}} \otimes \S_{m,n}} & 
p^{\otimes 0} \otimes p^{\otimes n} \dar[equals] \\
p^{\otimes 0} \otimes p^{\otimes m} \dar[equals,swap]{\wr} \ar[slash, ""{name=T, above}, ""{name=U, below}]{r}[description]{\S_{0,0} \otimes \S_{m,n}} & 
p^{\otimes 0} \otimes p^{\otimes n} \dar[equals,swap]{\wr} \\
p^{\otimes m} \rar[slash, ""{name=V, above},swap]{\S_{m,n}} & p^{\otimes n}
\arrow[equals,shorten=5,from=Q,to=R,swap,"\wr"]
\arrow[Rightarrow,shorten=5,from=S,to=T]
\arrow[Rightarrow,shorten=5,from=U,to=V]
\end{tikzcd}\!\!\!\!\!\!\! = \;\begin{tikzcd}[column sep=scriptsize]
p^{\otimes m} \dar[equals] \rar[slash,""{name=S, below}]{\S_{m,n}} & p^{\otimes n} \dar[equals] \\
p^{\otimes m} \rar[slash, ""{name=T, above},swap]{\S_{m,n}} & p^{\otimes n}
\arrow[equals,shorten=5,from=S,to=T]
\end{tikzcd}\; = \!\!\!\!\!\!\!\begin{tikzcd}[column sep={120,between origins}]
p^{\otimes m} \dar[equals,swap]{\wr} \rar[slash,""{name=Q, below}]{\S_{m,n}} & p^{\otimes n} \dar[equals,swap]{\wr} \\
p^{\otimes m} \otimes p^{\otimes 0} \dar[equals] \ar[slash, ""{name=S, below}, ""{name=R, above}]{r}[description]{\S_{m,n} \otimes \idcoalg{p^{\otimes 0}}} & 
p^{\otimes n} \otimes p^{\otimes 0} \dar[equals] \\
p^{\otimes m} \otimes p^{\otimes 0} \dar[equals,swap]{\wr} \ar[slash, ""{name=T, above}, ""{name=U, below}]{r}[description]{\S_{m,n} \otimes \S_{0,0}} & 
p^{\otimes n} \otimes p^{\otimes 0} \dar[equals,swap]{\wr} \\
p^{\otimes m} \rar[slash, ""{name=V, above},swap]{\S_{m,n}} & p^{\otimes n}
\arrow[equals,shorten=5,from=Q,to=R,swap,"\wr"]
\arrow[Rightarrow,shorten=5,from=S,to=T]
\arrow[Rightarrow,shorten=5,from=U,to=V]
\end{tikzcd}\hspace{-.7cm}
\end{equation}
\end{itemize}
\end{definition}

\section{Proofs of Dynamic Structure}\label{proofs}

We now proceed to prove that the coalgebras and structure maps defined above for organized predictions and gradient descent form dynamic structures. In each case, it suffices to show that the structure maps on states preserve coalgebra structure, and that the equations in \cref{operadequations} or \cref{PROequations}, respectively, are satisfied.

\begin{proof}[Proof of \cref{predictionadaptive}]
The operad equations are all satisfied as $\Delta^+$ is known to be an operad, and morphisms of coalgebras are entirely determined by a function between the state sets. It then remains only to show that the identitor and compositor as defined in \cref{sec.kelley} commute with actions and updates. This is clearly true for the identitor as it is an isomorphism, so we focus on the compositor.

For the compositor to commute with actions on positions is the claim that $\Delta^+_X$ is an algebra for the operad $\Delta^+$; it means that for $\mu \in \Delta^+_N$, $\nu_1 \in \Delta^+_{M_1}$, ..., $\nu_N \in \Delta^+_{M_N}$, and $\gamma_{i,j} \in \Delta^+_X$ for $i=1,...,N$ and $j=1,...,M_i$, we have
\[
\sum_i \mu_i \left(\sum_j \nu_j \gamma_{i,j}\right) = \sum_{i,j} (\mu_i \nu_{i,j})\gamma_{i,j},
\]
which is clearly the case.

The compositor commutes with actions on directions because in $(\S_{M_1} \otimes \cdots \otimes \S_{M_N}) \then \S_N$ the action of $(\nu_1,...,\nu_N,\mu)$ sends an outcome 
\[
x \in X = p_X[\mu \cdot (\nu \cdot \gamma)]
\]
to 
\[
(x,...,x) \in X^N = p_X^{\otimes N}[\nu_1 \cdot \gamma_1,...,\nu_N \cdot \gamma_N]
\]
and then to 
\[
(x,...,x) \in X^{\sum_i M_i} = p_X^{\otimes \sum_i M_i}[\gamma_{1,1},...,\gamma_{N,M_N}],
\]
while in $\S_{\sum_i M_i}$ the action of $\mu \circ \nu$ sends $x \in X$ to $(x,...,x) \in X^{\sum_i M_i}$ directly.

It then only remains to show that the compositor commutes with updates. Using the shorthand notation $\nu = (\nu_1,...,\nu_N)$ and $\gamma = (\gamma_1,...,\gamma_N) = (\gamma_{1,1},...,\gamma_{N,M_N})$ already employed above, to show that the composite of the updates of $\mu,\nu$ agrees with the update of the composite $\mu \circ \nu$ amounts to the equation
\begin{equation}\label{trustupdate}
\gamma(x) * (\mu \circ \nu) = \big((\nu \cdot \gamma)(x) * \mu\big) \circ \big(\gamma(x) * \nu\big)
\end{equation}
for any $x \in X$. Here $\nu \cdot \gamma$ denotes $(\nu_1 \cdot \gamma_1,...,\nu_N \cdot \gamma_N)$ and $\gamma(x) * \nu$ denotes $(\gamma_1(x) * \nu_1,...,\gamma_N(x) * \nu_N)$. On the $(i,j)$-component of these distributions, \eqref{trustupdate} unwinds to
\[
\frac{\gamma_{i,j}(x)(\mu_i\nu_{i,j})}{\sum_{i',j'}\gamma_{i',j'}(x)(\mu_{i'}\nu_{i',j'})} = \left(\frac{\sum_{j'}(\nu_{i,j'}\gamma_{i,j'}(x))\mu_i}{\sum_{i'}\sum_{j'}(\nu_{i',j'}\gamma_{i',j'}(x))\mu_{i'}}\right)\left(\frac{\gamma_{i,j}(x)\nu_{i,j}}{\sum_{j'}\nu_{i,j'}\gamma_{i,j'}(x)}\right),
\]
which is easily seen to hold by extracting $\mu_i$ from the first fraction on the right hand side and cancelling the sums over $j'$.
\end{proof}

\begin{proof}[Proof of \cref{gradientadaptive}]
The unit and associativity equations follow immediately from associativity and unitality of addition, cartesian products, and function composition. The interchange equations follow from the preservation of 0 under addition and identity functions under cartesian products, the analogous interchange property of function composition and cartesian products of functions, and the fact that the compositors and productors act the same way on the parameters and their dimension. 

It then remains only to show that the identitors, compositors, and productors are morphisms of coalgebras. This is immediate for the productors, as each component of the action and update functions respects the cartesian products of functions and parameters that define them, so we proceed only for the identitors and compositors.

For the identitors, the state $(0,\id_{\rr^n},0)$ acts as the identity function on $\rr^n$ and on directions by the transpose of its derivative, which is also the identity. The updates in the coalgebras $\S_{n,n}$ only modify the parameter $p$, so as the parameter here is 0-dimensional this state is never changed by the update function, as is the case in the coalgebra $\idcoalg{t^{\otimes n}}$. Therefore this function is a map of coalgebras.

The compositors preserve the component of the action on positions as, for states 
\[(L \in \nn, f\colon \rr^{L+\ell}, p \in \rr^L) \quad \textrm{and} \quad (M \in \nn, g\colon \rr^{M+m}, q \in \rr^M),\]
we have 
\[g(q,-) \circ f(p,-) = g(-,f(-,-))(q,p,-).\]
This may seem like a trivial rewriting, but it illustrates how the compositor was defined in order for the action to be preserved, as on the left we have the composite of the actions on positions as in $\S_{\ell,m}\then\S_{m,n}$, and on the right we apply the compositor and take the action of the resulting state in $\S_{\ell,n}$.

To show that the compositor preserves both the action on directions and the update we note that by the chain rule, for $x \in \rr^\ell$ and $z \in T_{g(q,f(p,x))}$,
\[D\left( g(-,f(-,-)) \right)^\top z = Df^\top \cdot \pi_m(Dg^\top \cdot z) \in T_{(p,x)}\rr^{L+\ell}.\]
Applying $\pi_\ell$ to both sides above shows that the compositor preserves the action on directions, as on the left we have the action on directions after applying the compositor and on the right we have the composition of the actions of $(L,f,p)$ and $(M,g,q)$ on directions as in $\S_{\ell,m}\then\S_{m,n}$.

Finally for updates, we observe by the chain rule that the update rule in $\S_{\ell,n}$ agrees with that in $\S_{\ell,m}\then\S_{m,n}$ under the compositor, as either way for $x,z$ as above the composite state of $(L,f,p)$ and $(M,g,q)$ updates to 
\[\left( M+L,g(-,f(-,-)),\left( p + \epsilon \pi_L (Df^\top \cdot \pi_m(Dg^\top \cdot z)),q + \epsilon \pi_M(Dg^\top \cdot z) \right) \right).\qedhere\]
\end{proof}

\end{document}